\input Tex-document.sty
\input psfig.sty
\def\R{I\!\! R}
\def\Z{\it Z \!\!\! Z}

\def\firstpage{279}

\pageno=279

\def\shorttitle{Action Minimizing Solutions of the Newtonian $n$-body
Problem}
\def\shortauthor{A. Chenciner}
\font\boldmathHuge = cmmib10 scaled 1728

\title{\centerline{Action Minimizing Solutions of the } 
\centerline{Newtonian {\boldmathHuge n}-body Problem: From } \centerline{Homology to Symmetry}}

\author{A. Chenciner\footnote{\eightrm $^\ast$}{\eightrm Astronomie et
Syst\`emes Dynamiques, IMCCE, UMR 8028 du CNRS, 77, avenue Denfert-Rochereau,
75014 Paris, France \& D\'epartement de Math\'ematiques,
Universit\'e Paris VII-Denis Diderot, 16, rue Clisson, 75013 Paris, France.
E-mail: chencine@bdl.fr}}

\vskip -5mm

\author{({\it A la m\'emoire de Nicole Desolneux})}

\vskip 7mm

\centerline{\boldnormal Abstract}

\vskip 4.5mm

{\narrower \ninepoint \smallskip
An action minimizing path between two given configurations, spatial or planar,
of the $n$-body problem is always a true -- collision-free -- solution.
Based on a remarkable idea of Christian Marchal, this theorem implies the existence of new
``simple" symmetric periodic solutions, among which the Eight for 3 bodies, the Hip-Hop for 4 bodies and their
generalizations.\smallskip

\vskip 4.5mm

\smallskip

\smallskip }

\vskip 10mm

\head{0. Introduction}

Finding periodic geodesics on a riemannian manifold as length minimizers in a fixed non-trivial homology or homotopy
class is commonplace lore. Advocated by Poincar\'e [P] as early as 1896, the search for periodic solutions of a
given period $T$ of the $n$-body problem as action minimizers in a fixed non-trivial homology or homotopy
class is rendered difficult by the possible existence of collisions due to the relative weakness of the newtonian
potential: the action of a solution stays finite even when some of the bodies are colliding. Very few results are
available: among them Gordon's characterization of Kepler solutions [G] for 2 bodies in $\R^2$, Venturelli's
characterization of Lagrange equilateral solutions [V1] for 3 bodies in
$\R^3$, Arioli, Gazzola and Terracini's characterization  of retrograde
Hill's orbits [AGT] for the restricted 3-body problem in $\R^2$. In particular,
no truly new solution of the
$n$-body problem was found in this way; indeed, these results confirm the view that the action-minimizing
periodic solutions are the ``simplest" ones in their class.
\smallskip
The action minimization method has recently been given a new impetus by
the replacement of the topological constraints by symmetry ones. This
idea was first introduced by the italian school [C-Z][DGM][SeT] as
another mean of  forcing coercivity of the problem, i.e. forbidding a
minimizer to be ``at infinity". The bodies were forced to occupy, after
half a period, a position symmetrical of the original one with respect
to the center of mass of the system. It is proved in [CD] that in a
space of even dimension, say $\R^2$, the minimizers in this symmetry
class include relative equilibrium solutions  (i.e. solutions which are
``rigid body like"); moreover all minimizers are of this form provided a
certain ``finiteness" hypothesis is verified (see [C3]). Such relative
equilibria can occur only for the so called {\it central configurations}
[AC], the most famous of which is Lagrange equilateral triangle.

Recently, a new type of symmetry was considered, which originates in the invariance of the Lagrangian under
permutations of equal masses. This has led to the discovery of a whole world of
new solutions in the case when all the bodies have
the same mass. The most surprising ones are the ``choreographies" whose name, given by Carles Sim\'o, fits the beautiful figures they display
on the screen in animated computer experiments ([CGMS],[S2]). Referring to my survey
article [C3] for a bibliography and a description of the few cases in which
existence proofs are available (the Hip-Hop [CV] for 4 bo\-dies in
$\R^3$, the Eight [CM] for 3 bodies in
$\R^2$, Chen's solutions [Ch] for 4 bodies in
$\R^2$),  I mainly address here a powerful theorem which
solves completely the collision problem for the fixed ends problem in the case of
 arbitrary masses. This is pertinent because, as we shall see, it allows one to prove the
existence of collision-free minimizers under well chosen symmetry constraints. This
theorem is the result of the efforts of Richard Montgomery, Susanna Terracini, Andrea
Venturelli [V2],   and, for the last -- fundamental
-- stone, Christian Marchal [M2] [M3]. I present here a complete proof and, in particular, a simplified
version of  Marchal's remarkable idea,
which avoids numerical computations. I discuss also  new applications to minimization
under symmetry constraints and open problems.

\noindent {\bf Notations.} By a {\it configuration} of $n$ bodies in an
euclidean space $(E,\left< \right>)$ we understand an $n$-tuple $x=(\vec
r_1,\vec r_2,\ldots \vec r_n)\in E^n$. The {\it configuration space} of
the $n$-body problem is the quotient of the set of configurations by the
action of translations (see [AC]). It may be identified as in [C3] with
the set ${\cal X}$ of configurations whose center of mass $\vec
r_G=(\sum_{i=1}^n{m_i})^{-1}\sum_{i=1}^n{m_i\vec r_i}$ is at the origin.
It is endowed with the ``mass scalar product" $(\vec r_1,\ldots,\vec
r_n)\cdot (\vec s_1,\ldots,\vec s_n)=\sum_{i=1}^n{m_i\left<(\vec
r_i-\vec r_G),(\vec s_i-\vec s_G )\right>}$. The {\it non-collision}
configurations -- the ones such that no two  bodies $\vec r _i$ coincide
-- form an open dense subset $\hat{\cal X}$ of ${\cal X}$. The functions
$I=x\cdot x, \; J=x\cdot y,\; K=y\cdot y$, defined on the {\it phase
space} $\hat{\cal X}\times {\cal X}$ (whose elements are noted $(x,y)$)
are the basic isometry-invariants of the $n$-body problem They are
respectively the {\it moment of inertia} of the configuration with
respect to its center of mass, half its time derivative and twice the
{\it kinetic energy} in a galilean frame which fixes the center of mass.
The {\it potential function} (opposite of the potential energy), the
{\it Hamiltonian} (=total energy) and the {\it Lagrangian} are
respectively defined by
$$U=\sum_{i<j}{m_im_j||\vec r_i-\vec r_j||^{-1}},\; H={1\over 2}K-U, \; L={1\over 2}K+U.$$
In terms of the gradient $\nabla$ for the mass metric, the equations of
the $n$-body problem,
$$m_i\ddot r_i(t)=\sum_{j\not=i}m_im_j{\vec r_j(t)-\vec r_i(t)\over \|\vec r_j(t)-\vec r_i(t)\|^3},\quad
i=1,\ldots,n,$$  can be written $\, \ddot x=\nabla U(x)$. They
are the  Euler-Lagrange equations of the action, which to a path $x(t)$ associates the real number
$${\cal A}_T\bigl(x(t)\bigr)=\int_0^T{L\bigl(x(t),\dot x(t)\bigr)dt}.$$
\noindent {\bf Remark.} In the perturbations, we shall not
bother about fixing the center of mass because replacing $K=\sum_{i=1}^n{m_i||\vec v_i-\vec
v_G||^2}$ by
$\sum{m_i||\vec v_i||^2}$ only increases the action.

\head{1. The fixed-ends problem}

\proclaim Question. {\rm Given two configurations, -- possibly
with collisions -- of $n$ point masses in $\R^3$ (resp. $\R^2$)
and a positive real number $T$, does there exist a solution of the
Newtonian $n$-body problem which connects them in the time $T$ ?}

A natural way of looking for a solution is to seek for a minimizer of  the action
${\cal A}_T(x)$ over the space
$\Lambda_0^T(x_i,x_f)$ of paths $x(t)$ in the configuration space $\hat{\cal X}$ which start at
time
$0$ in the configuration
$x_i$ and end at time
$T$ in the configuration $x_f$. For the integral to be defined, it is natural to work in the
Sobolev space of paths which are square integrable together with their
first derivative in the sense of distributions.

The main problem, already mentioned by  Poincar\'e in 1896 (see [P]
where he introduces the method in a slightly different context), is that
a minimum could  well be such that, for a non-empty set of instants
(necessarily of measure zero), the system undergoes a collision of two
or more bodies, which prevents it form being a true solution (see [C3];
existence is not a problem here because fixing the ends gives
coercivity). At an isolated collision time, the renormalized
configuration is known to be approaching the set of central
configurations (the ones which admit homothetic motion [C2]) but very
little is understood of these configurations for more than 3 bodies.
Continuous families of such configurations could exist (the ``finiteness
problem") and even if they didn't, there would be no garantee that at
collision the renormalized configuration has a limit : it might have one
only modulo rotations ( the ``infinite spin problem"). Nevertheless, we
prove the \proclaim Theorem. {\it A minimizer of the action in
$\Lambda_0^T(x_i,x_f)$ is collision-free on the whole open interval
$]0,T[$. Hence, the answer to the Question is yes, both in $\R^3$ and
$\R^2$.}

In the next paragraph, Marchal's idea to prove that isolated collisions
do not occur in a (local) minimizer is explained on the Kepler problem.
If the finiteness problem is supposed to be solved, it works in the same
way for the general $n$-body problem (surprisingly, the infinite spin
problem is irrelevant). We then address the finiteness problem with
Terracini's technique of {\it blow up}, which reduces the problem of
isolated collisions to the case of parabolic homothetic solutions;
finally we show, following Montgomery and Venturelli, that accumulation
of collisions do not occur in a minimizer provided no subclusters
collide. The theorem then follows by induction on the number of bodies
involved in a collision.
\smallskip\noindent {\bf Remark.} A similar assertion, based on numerical experiments, was
made by Tiancheng Ouyang in Guanajuato (Hamsys, march 2001) but no proof
has yet appeared.

\head{2. The Kepler problem as a model for the study of isolated collisions}

The case of two bodies contains already many ingredients of the general situ\-ation. As is well-known, the 2-body
problem is equivalent to the problem of a particle attracted to a fixed center 0, the so-called Kepler problem (or
1-fixed center problem). We call {\it collision-ejection} a solution in which the particle follows a straight line
segment from its initial position $\vec r_i$ to the attracting center and (possibly) another straight line segment
from the attracting center to its final position $\vec r_f$.

\proclaim A test assertion. {\it A collision-ejection solution of
the Kepler problem does not minimize the action in the Sobolev
space $\Lambda_0^T(\vec r_i,\vec r_f)$ of paths joining $\vec
r_i$ to $\vec r_f$.}

At least four proofs may be given of the truth of this assertion but only the fourth one using
Marchal's idea is robust enough to lead to complete generalization.
In the first one, we use the explicit knowledge of the solutions of the 2-body problem [A1] to identify the
minimizers with the ``direct" arcs of solution, not going ``around" the attracting center (this arc is uniquely
determined provided
$\vec r_i$, O  and $\vec r_f$ do not lie on a line in this order). In the second one,  we
find a ``simple" path without collision (straight line, circle, uniform
motion) which has lower action. In the third one, supposing that a minimizer $\vec r(t)$ has a collision with
the fixed center at  time $0$, we
find a local deformation $\vec r_\epsilon(t)=\vec r(t)+\epsilon\varphi(t)\vec s$, which has
lower action and no collision. Such deformations were used by
many people, including Susanna Terracini, Gianfausto Dell'Antonio, Richard Montgomery and Christian Marchal.
If we chose, with Montgomery,
$\varphi(t)=1$ if $0\leq t\leq \epsilon^{3\over 2},\;
\varphi(t)=\epsilon^{-1}(\epsilon^{3\over 2}+\epsilon-t)$ if $\epsilon^{3\over 2}\leq t\leq \epsilon^{3\over
2}+\epsilon\;$ and $\varphi(t)=0$ if $t\geq
\epsilon^{3\over 2}+\epsilon$,
the gain in action is
$c\sqrt{\epsilon}\left(1+O(\sqrt{\epsilon}\log(1/\sqrt{\epsilon}))\right)$ provided the unit vector $\vec s$ is
well chosen. We come to the fourth proof, for which we must distinguish two cases according to the dimension
of the ambient space.
\goodbreak
\smallskip
\noindent {\bf (i) The case of $\R^3$.} Let $t\mapsto \vec r(t)$ be a collision-ejection solution of the
Kepler problem,
$\ddot{\vec r}(t)=-\vec r(t)/|\vec r(t)|^3$, such that $\vec r(-T')=\vec r_i,\,\vec r(0)=0,\,\vec r(T)=\vec r_f,\,
T,T'>0$.
We consider the following family of continuous deformations of $\vec r(t)$, parametrized by an element $\vec
s$ of the unit sphere $S^2$ in $\R^3$ : if $R'(t)=(1+{t\over T'})\rho$ and $R(t)=(1-{t\over T})\rho$,
$$\vec r_{\vec s}(t)=\vec r(t)+R'(t)\vec s\;\;\hbox{if}\;\; -T'\leq t\leq 0,\quad
\vec r_{\vec s}(t)=\vec r(t)+R(t)\vec s\;\;\hbox{if}\;\; 0\leq t\leq T.$$
It is a simplification of Marchal's original choice but the idea is the same : to show that the action ${\cal A}$ of
$\vec r(t)$ is strictly bigger than the average ${\cal A}_m=\int_{S^2}{{\cal A}(\vec r_{\vec s}(t))d\sigma}$,
where $d\sigma$ denotes the normalized area form, that is the unique rotation invariant probability measure
on
$S^2$. This will imply the existence of at least one direction $\vec s$ for which $\vec r_{\vec s}(t)$ has lower
action than
$\vec r(t)$ (because the set of good $\vec s$ has positive measure we could choose $\vec s$ so that
$\vec r_{\vec s}$ is collision-free but this is irrelevant).

The linearity of the integral and the similar behaviour of  ejection
and collision allow {\it to replace in the proof $\vec r(t)$ and $\vec r_{\vec s}(t)$ by
their restrictions to the interval $[0,T]$}.
Moreover, it follows from the ``blow-up" method (see 3.2) that it is enough to consider a {\it
parabolic} solution $\vec r(t)$, that is
$\vec r(t)=\gamma t^{2\over 3}\vec c$, with $\gamma=(9/2)^{1\over 3}$ if $|\vec c|=1$.
\smallskip
By Fubini theorem applied to the positive integrand,
$${\cal A}_m=\int_{S^2}d\sigma\int_0^T\left({|\dot{\vec r_{\vec s}}|^2\over
2}+{1\over |\vec r_{\vec s}|}\right)dt=\int_0^Tdt\int_{S^2}\left({|\dot{\vec r_{\vec s}}|^2\over
2}+{1\over |\vec r_{\vec s}|}\right)d\sigma,\quad\hbox{ and }$$
$${\cal A}_m-{\cal A}=\int_0^T{dt\left[{\dot R(t)^2\over 2}+\int_{S^2}{\dot R(t)\vec s\cdot\vec r(t)\,  d\sigma}\right]}
+\int_0^T{dt\left[\int_{S^2}{d\sigma\over |\vec r_{\vec s}(t)|}-{1\over |\vec r(t)|}\right]}.$$
The first integral reduces to $\, {1\over 2}\int_0^T{\dot R(t)^2}dt=\rho^2/2T$ because of the antisymmetry in
$\vec s$ of the scalar product. The second is the
difference in potential resulting from the replacement of the particle $\vec r(t)$ by a
homogeneous hollow sphere of the same mass and increasing radius $R(t)$.
Because of the harmonicity of Newton potential in $\R^3$, the potential $U_0(\vec r, R):=\int{d\sigma\over |\vec r-R\vec
s|}=\int{d\sigma\over |\vec r+R\vec s|}$ of a homogeneous hollow sphere of radius $R$  is
$$U_0(\vec r,R)={1\over R}\quad \hbox{if}\quad |\vec r|\leq
R,\quad U_0(\vec r,R)={1\over |\vec r|}\quad \hbox{if}\quad |\vec r|\geq
R.$$ If $0$ leaves this sphere at time $t_0$, $|\vec r(t_0)|=R(t_0)$,
i.e. $\rho=\gamma t_0^{2\over 3}+O(t_0^{5\over 3})$, and
$${\cal A}-{\cal A}_m={\rho^2\over 2T}+\int_0^{t_0}{\left[{1\over R(t)}-{1\over |\vec
r(t)|}\right]dt}=-{2\over \gamma}t_0^{1\over 3}+O(t_0^{4\over 3 })\leq 0\;\hbox{if $\rho$, hence $t_0$, is small.
}$$

\noindent {\bf (ii) The case of $\R^2$.} The Newtonian potential is not
harmonic in $\R^2$ and this makes things somewhat more complicated.
Marchal proposes to replace the sphere by a disk of radius $R$ endowed
with the projection $\sigma(\theta,x)=1/\bigl(2\pi
R\sqrt{R^2-x^2}\bigr)$ (in polar coordinates) of the uniform density on
the sphere of the same radius. The potential fonction $U_0(\vec r,R)$ of
such a disk (total mass 1) may be recovered from the general computation
done, via complex function theory, for a thin elliptic plate with a
given density which is constant on homothetic ellipses (see [B]):
$$U_0(\vec r,R)={\pi\over 2R}\quad \hbox{if}\quad |\vec r|\leq
R,\quad
U_0(\vec r,R)={1\over R}\arcsin\bigl({R\over |\vec r|}\bigr)\quad \hbox{if}\quad |\vec r|\geq
R.$$

It does not coincide any more, but asymptotically, with Newton's
potential $1/|\vec r|$ of the center of mass outside the disk but it is
still constant in the interior and the proof works as well as in the
spatial case:  as $\arcsin(x)\leq x+({\pi\over 2}-1)x^3$ for $x\geq 0$,
the difference in actions between the mean of the modified actions when
$\vec s$ belongs to the unit disk and the original becomes
$$\eqalign{{\cal A}_m-{\cal A}&={\rho^2\over 2T}+\int_0^{t_0}{\left[{\pi\over 2R(t)}-{1\over |\vec
r(t)|}\right]dt}+\int_{t_0}^T{\left[{1\over R(t)}\arcsin\bigl({R(t)\over |\vec r(t)|}\bigr)-{1\over |\vec
r(t)|}\right]dt}\cr
&\leq {\rho^2\over 2T}+\left[-{\pi T\over 2\rho}\log(1-{t\over T})-{3\over \gamma}t^{1\over
3}\right]_0^{t_0}+\int_0^{t_0}{({\pi\over 2}-1)\rho^2(1-{t\over T})^2{1\over \gamma^3t^2} dt}\cr
&={\gamma^2t_0^{4\over 3}\over 2T}\bigl(1+O(t_0)\bigr)+({\pi\over 2}-3){1\over
\gamma}t_0^{1\over 3}\bigl(1+O(t_0)\bigr)+({\pi\over 2}-1){1\over \gamma}t_0^{1\over
3}+O\bigl(t_0^{4\over 3}\log({1\over t_0})\bigr)\cr &=(\pi -4){1\over
\gamma}t_0^{1\over 3}+O\bigl(t_0^{4\over
3}\log({1\over t_0})\bigr)\leq 0\;\hbox{ for $\rho$, hence $t_0$, small.}}$$

\head{3. Proof of the theorem}

\noindent {\bf 3.1 The induction.} We define the following statements about a minimizer $x(t)$:
\smallskip
$(I_p)$ If a collision of $p$ bodies occurs in $x(t)$ for $t\in]0,T[$, it is isolated.

$(II_p)$ No collision of $m\leq p$ bodies occurs in $x(t)$ for $t\in]0,T[$.
\smallskip
\noindent In 3.2 we  prove that $(I_p)$ implies that no collision of $p$ bodies occurs in $]0,T[$, hence that
$(II_p)$ and $(I_{p+1})$ imply $(II_{p+1})$.
In 3.3  we prove that
$(II_p)$ implies $(I_{p+1})$.
As $(II_1)$ is empty, it implies $(I_2)$, hence $(II_2)$, etc... up to  $(II_n)$ which
is the conclusion. If $k$-collisions are present in $x_i$ or $x_f$ but not $j$-collisions for $j<k$, the
induction proves that $j<k$-collisions are absent. The next step proves that $k$-collisions,
including the ones at the ends, are isolated and everything goes through.
\smallskip
\noindent {\bf Remark.} The induction may succeed because a $p$-body collision cannot be a
limit of $q$-body collisions with $q>p$. Still, accumulation of collisions involving bodies in different clusters
could {\it a priori} occur, e.g. a sequence of double collisions 23, 12, 34, 23, 13, 24, 23, ... converging to a
quadruple collision 1234, or even a converging sequence of such sequences. Induction on the number of
bodies in the clusters fortunately avoids having to deal with such problems.

\medskip

\noindent {\bf 3.2 Elimination of isolated collisions.}
\smallskip
\noindent {\bf 3.2.1 The blow-up technique.} This technique was
introduced by S. Terracini and developped in the thesis of A.
Venturelli [V2]. It is based on the homogeneity of the potential
(compare [C2]). It allows proving the
\proclaim Proposition. {\it
If a minimizer $x(t)$ of the fixed ends problem for $n$-bodies
possesses an isolated collision of $p\leq n$ bodies, there is a
parabolic  (i.e. zero energy) homothetic collision-ejection
solution $\overline x(t)$ of the $p$-body problem which is also a
minimizer of the fixed ends problem.}

\noindent {\bf Proof.} To keep the exposition as simple as possible, I describe the case of a
total collision. In the general case of partial (and possibly simultaneous) collisions, everything goes
through in the same way because the blow up sends all bodies not concerned by the collision to infinity (for
more details, see [V2]).
\smallskip
Assuming that the collision occurs at $t=0$, we define
$x^\lambda(t)=\lambda^{-{2\over 3}}x(\lambda t)$ for $\lambda>0$.
If $x(t)$ is a solution of the $n$-body problem, so is $x^\lambda(t)$. Moreover, for any path $x(t)$ in
$\Lambda_{T_1}^{T_2}(x_i,x_f))$, the path $x^\lambda(t)$ belongs to
$\Lambda_{T_1}^{T_2}(x^\lambda(T_1),x^\lambda(T_2))$ and its action is equal to $\lambda^{-{1\over 3}}$
times the action of the restriction of $x(t)$ to the interval $[\lambda T_1,\lambda T_2]$. Hence, if
$x(t)$ is action minimizer in $\Lambda_{T_1}^{T_2}(x_i,x_f)$, so is $x^\lambda$ in
$\Lambda_{T_1}^{T_2}(x^\lambda(T_1),x^\lambda(T_2))$. Now, Sundman's estimates recalled above
imply that,
$\{x^\lambda, 0<\lambda<\lambda_0\}$ is bounded in $H^1([T_1,T_2],{\cal X})$, hence weakly compact, so
that there exists a sequence $\lambda_n\to 0$  such that
$x^{\lambda_n}$ converges weakly (and hence uniformly) in $H^1([0,T],{\cal X})$ to a solution
$\overline x$. One shows that $\overline x$ is made of a parabolic homothetic collision
solution followed by a parabolic homothetic ejection solution (the two central configurations involved are a
priori distinct). Moreover, it follows from the weak lower semi-continuity of the action that $\overline x$ is a
minimizer in
$\Lambda_{T_1}^{T_2}(\overline x(T_1),\overline x(T_2))$ (see [V2]).
\medskip
\noindent {\bf 3.2.2 The mean perturbed action.} We shall deal only with the
case of $\R^3$ and refer the reader to the Kepler case for the modifications needed in the case of $\R^2$.
Thanks to ``blow up", we may suppose that our minimizer $x(t)$ is a parabolic homothetic collision-ejection
solution
$x(t)=(\vec r_1(t),\cdots,\vec r_p(t))=x_0|t|^{2\over 3}$ of the
$p$-body problem. As in the Kepler case, we may restrict to the ejection part, corresponding to $t\in [0,T]$.
One studies deformations of $x(t)$ of the form
$$x_{\vec s}^k(t)=\bigl(\vec r_1(t),\ldots,\vec r_k(t)+R(t)\vec
s,\ldots,\vec r_p(t)\bigr),$$
where, as before, $R(t)=(1-{t\over T})\rho$ with
$\rho$ a small positive real number and
$\vec s$ belongs to the unit sphere.
The same computation as in the Kepler case leads to an average action ${\cal A}_m^k$ such that
$${\cal A}_m^k-{\cal A}\leq {m_k\over 2}{\rho^2\over T}+
\sum_{j\not=k,\, j\leq p}{m_jm_k\int_0^{t_{jk}}\left[{1\over R(t)}-{1\over r_{jk}(t)}\right]dt},$$
where $r_{jk}=|\vec r_k-\vec r_j|$ and $t_{jk}$ is defined by $r_{jk}(t)=R(t)$
(the inequality sign comes from the fact that the deformations do not keep the center of mass fixed).

As $r_{jk}(t)=c_{jk}t^{2\over 3}$, one concludes as in the Kepler case that  ${\cal A}_m^k-{\cal A}<0$.
\smallskip
\noindent {\bf Remark.} We could have dispensed with ``blow up" in case similitude classes of central
configurations were isolated but certainly not otherwise. This is because, the best control Sundman's theory
may give us on the asymptotic behaviour of the colliding bodies is that their moment of inertia $I_c$ with
respect to their center of mass and their potential $U_c$ are respectively
equivalent to $I_0t^{4\over 3}$ and $U_0t^{-{2\over 3}}$ (see [C2]).
This implies the existence, for
$1\leq j<k\leq p$,  of $0<a_{jk}\leq
b_{jk}$ such that for $t$ small enough, one has
$a_{jk}t^{2\over 3}\leq r_{jk}(t)\leq b_{jk}t^{2\over 3}$.
It follows that
$${\cal A}_m^k-{\cal A}\leq {m_k\over 2T}b_{jk}^2t^{4\over 3}+O(t_{jk}^{7\over 3})-\sum_{j\not=k,\,
j\leq p}m_jm_k\left(\left[-{1\over a_{jk}}+{3\over b_{jk}}\right]t_{jk}^{1\over 3}+o(t_{jk}^{1\over
3})\right).$$  If similitude classes of central configurations are isolated, there is a limit shape and we may take
$a_{jk}$ and $b_{jk}$ as close as we wish. Otherwise we cannot conclude.
\medskip

\noindent {\bf 3.3 The elimination of non-isolated collisions.}
It remains to prove that $(II_p)$ implies $(I_{p+1})$. We use
energy considerations, an idea which goes back to R. Montgomery
and was further developed in Venturelli's thesis [V2].

\proclaim Proposition. {\it Let $x(t)$ be a minimizer of the fixed
ends problem. If $x(t)$ has no $p$-body collisions for $p<p_0$,
collisions of $p_0$ bodies are isolated.}

\smallskip
\noindent {\bf Sketch of proof.} I shall give the proof in the case of a total collision (i.e. $p_0=n$) and
then explain what has to be changed in the general case.
\smallskip
(i) Using the behavior of the action under reparametrization, let us prove that the energy stays constant along a
minimizer, whatever be the collisions. For this let us consider variations
$x_\epsilon(t)$ of the form $x_\epsilon(t)=x(\varphi_\epsilon(t))$ where $t\mapsto \tau=\varphi_\epsilon(t)$
is a differentiable family of diffeomorphisms of $[0,T]$ starting from $\varphi_0(t)\equiv t$ :
  $${\cal A}(x_\epsilon)=\int_0^T\left({||\dot x_\epsilon(t)||^2\over
2}+U(x_\epsilon(t))\right)dt =\int_0^T\left({1\over \lambda_\epsilon(\tau)}{||\dot x(\tau)||^2\over
2}+\lambda_\epsilon(\tau)U(x(\tau))\right)d\tau,$$ where
$\lambda_\epsilon=dt/d\tau=1/\dot\varphi_\epsilon(\varphi_\epsilon^{-1}(\tau))$.  The derivative at
$\epsilon=0$ of
$a(\epsilon)={\cal A}(x_\epsilon)$ is
$${da\over d\epsilon}(0)=\int_0^T\left({||\dot x(\tau)||^2\over
2}-U(x(\tau))\right)\delta\lambda(\tau)\, d\tau=\int_0^T{H\bigl(x(\tau),\dot
x(\tau)\bigr)\delta\lambda(\tau)d\tau},$$
where $\delta\lambda(\tau)={d\lambda_\epsilon(\tau)\over
d\epsilon}|_{\epsilon=0}$. As the variations $\delta\lambda$ satisfy the constraint
$\int_0^T{\delta\lambda(\tau)d\tau}=0$, which comes from the fact that
$\int_0^T{\lambda_\epsilon(\tau)d\tau}=T$, we get that there exists a real constant $c$ such
$H\bigl(x(\tau),\dot x(\tau)\bigr)=c$ wherever it is defined.

\smallskip
(ii) Let $t_0$ be an instant at which total collisions accumulate. Let us chose two sequences $(a_n)$ and $(b_n)$
of instants of total collision which converge to
$t_0$ and are such that no total collision occurs in the open intervals $]a_n,b_n[$. The moment of inertia $I$ of
the system with respect to its center of mass is equal to zero at each of the instants $a_n$ or $b_n$ and hence has
at least one maximum $\xi_n$ in the interval $]a_n,b_n[$. As no partial collision occurs,  the motion is regular in
each of these intervals and at each such maximum, the second time-derivative
$\ddot I(\xi_n)$ has to be non positive. But the value
$U(\xi_n)$ of the potential function tends to $+\infty$ as $n\to +\infty$, while the energy $H$
stays constant. One then deduces from the Lagrange -Jacobi relation
$\ddot I=4H+2U$
that $\ddot I(\xi_n)\to +\infty$, which is a
contradiction.
\smallskip
In the general case, when $\mu$ is some cluster not containing all the bodies, the energy
$H_\mu$ of $\mu$ is no more constant but one can get from a refinement of the same proof that it is still an
absolutely continuous function of time as long as no collision occurs between a body of the cluster and a body of
the complementary cluster (see [V2]). This implies that $H_\mu$ stays locally bounded and allows the argument
of (ii) to work because, by hypothesis, no partial collision occurs in the cluster.

\head{4. Periodic solutions}

\noindent {\bf 4.1 Homological or homotopical constraints.} Going back
to the 1896 Note of Poincar\'e already alluded to, the idea of
constructing periodic solutions of the $n$-body problem as the
``simplest" (action minimizing) ones in a given homology or homotopy
class of the configuration space is very natural if one compares to the
construction of periodic geodesics as minimizing the length in a non
trivial homology or homotopy class. As already noticed by Poincar\'e,
this works beautifully in the so-called ``strong force problem",
corresponding to a potential in $1/r^2$ or stronger, where each
collision path has infinite action [CGMS].  Unfortunately, in the
Newtonian case, most of the time minimizers have collisions and hence
are not true periodic solutions [M]. This is already true in the planar
Kepler problem: it follows from Gordon's work [G] (see also [C3]) that
the only minimizers of the action among the loops of a fixed period $T$
whose index in the punctured plane is different from $0,\pm 1$, is an
ejection-collision one ! (for an analogue result in the planar
three-body problem, see [V1]).
\smallskip
In such cases, solving the fixed ends problem is of no use.
Among the cases where minimizers in a fixed homology or homotopy class have no collision are

1) Gordon's theorem for the planar Kepler problem when one fixes
the index to $\pm 1$ (resp. when one insists only on the index
being different from 0): a minimizer is any elliptic solution of
the given period.

2) The generalization [V1],[ZZ1] of Gordon's theorem to the planar
three-body problem whith homology class fixed in such a way that along a
period, each side of the triangle makes exactly one complete turn in the
same direction: a minimizer is any elliptic homographic motion
 of the equilateral triagle, of the given period.
\medskip
\noindent {\bf 4.2 Symmetry constraints.} In order to find ``new"
solutions as action minimizers, another type of constraints on
the loops must be introduced, which somewhat allows using fixed
ends type results. We ask the loops to be invariant under the
action of a finite group $G$. An invariant loop is completely
defined by its restriction to an interval of time on which $G$
induces no constraint. The restriction to such a``fundamental
domain" of a minimizer among $G$-invariant loops is a minimizer
of the fixed ends problems between its extremities. This leads to
a new collision problem: a minimizer could well have a collision
at the initial or final instant.
\smallskip
\noindent {\bf (i) Choreographies.} We first show, following
Andrea Venturelli, the

\proclaim Theorem. {\it A minimizer among $n$-choreographies has
no collision.}

Recall that the choreographies are fixed loops under the action of the
group ${\Z}/n{\Z}$ whose generator cyclically permutes equal-mass bodies
after one $n$-th of the period (see [CGMS]); hence a fundamental domain
can be chosen as any time interval of length $T/n$. If there were
collisions at the ends, one would get a contradiction with the theorem
by just shifting the fundamental domain to the right or to the left. One
can prove (using [CD]) that the regular $n$-gon is action minimizing
when it minimizes $\tilde U=I^{1\over 2}U$. But this is no more true for
$n\geq 6$. So, what is the min ?
\smallskip
\noindent {\bf (ii) Generalized Hip-Hops.} This works also for the ``italian" (anti)symmetry:

\proclaim Theorem.  {\it A minimizer among loops $x(t)$ in $\R^3$
such that $x(t+T/2)=-x(t)$ has no collision. Moreover, it is
never a planar solution.}

The last assertion comes from the
fact that a relative equilibrium $x(t)$ whose configuration $x_0$ minimizes $\tilde U=I^{1\over 2}U$ is always
a minimizer among the planar (anti)symmetric loops ([CD] and [C3]). But, applied to a variation
$z(t)=z_0\cos{2\pi t\over T}$ normal to the plane of $x(t)$, the Hessian of the action is easily seen [C4] to be
$$d^2{\cal A}(x(t))(z(t,z(t))=I_0^{-{1\over 2}}d^2\tilde U(x_0)(z_0,z_0)\int_0^T{\cos^2{2\pi t\over T}dt},$$
where $I_0=x_0\cdot x_0$. Now,  results of Pacella and Moeckel [Mo1] say
that one can always choose $z_0$ such that  $d^2\tilde
U(x_0)(z_0,z_0)<0$. Hence, a relative equilibrium ceases being a
minimizer in $\R^3$. This ends the proof because other possible
minimizers of the planar problem would have the same action as a
relative equilibrium (thanks to A. Venturelli for this remark). In
reference to  [CV], I propose to call {\it generalized Hip-Hops} these
minimizers. They are the best approximations I can think of in $\R^3$ to
the non-existing relative equilibria of non-planar central
configurations (recall that, according to [AC] such relative equilibria
exist in $\R^4$).
\smallskip
\noindent {\bf (iii) Eights with less symmetry.} As another example, we
prove the existence, for three equal masses, of solutions ``of the Eight
type" but with less a priori symmetry than the full symmetry group
$D_6=\{s,\sigma|s^6=1,\sigma^2=1,s\sigma=\sigma s^{-1}\}$ (see [C3]) of
the space of oriented triangles (``shape sphere" in [CM]). We consider
the subgroups $\Z/6\Z=\{s\}$ and $D_3=\{s^2,\sigma\}$.

\proclaim Theorem. {\it A minimizer among $\Z/6\Z$-invariant loops
has no collision. The same is true for a minimizer among
$D_3$-invariant loops.}

Instead of minimizing the action over one twelfth of the period between an Euler configuration at time 0 and
an isosceles one at time $T/12$ (see [CM]), one minimizes only over one sixth of the period: in the first case
from an isosceles configuration at time $t_0$ to a symmetric one at time $t_0+T/6$, in the second one from an
Euler configuration at time $0$ to another one at time $T/6$. Venturelli's trick of translating the fundamental
domain works in the first case where $t_0$ is arbitrary (a translation of time transforms a minimizer into a
minimizer) but not in the second one where, as for the initial
$D_6$-action, an Euler configuration can only occur at times which are integer multiples of $T/6$.
To prove the absence of collisions at the initial and the final instant in the second case, we notice that such
a collision is necessarily a triple (i.e. total) collision. If this happens, the action of the path  is greater than the
one of a homothetic ejection solution of equilateral type, a path which is not of the required type, but this is
irrelevant here. The conclusion follows because the action of this last path is itself greater than the one  of one
sixth of the ``equipotential model" (see [C3],[CM]).
If a minimizer among $\Z/6\Z$ or $D_3$ symmetric
loops possesses the whole
$D_6$ symmetry of the Eight is unknown.
\smallskip
\noindent {\bf (iv) The $P_{12}$ family.} Marchal discovered the
$P_{12}$ family, which continues the Eight solution in three-space up to
Lagrange equilateral solution, through choreographies in a rotating
frame [M1]. It is parametrized by an angle $u$ between $0$ and $\pi\over
6$:  the solution labeled by $u$ is supposed to minimize the action in
fixed time $T/12$ between configurations which are symmetric with
respect to a line $\Delta$ through the origin which contains body $0$
and configurations which are symmetric with respect to a plane $P$
through the origin which contains body 2 and makes angle $u$ with
$\Delta$.  We shall think of $\Delta$ as being horizontal and of $P$ as
being vertical (Figure 1). \midinsert

\centerline{\hbox{\psfig{figure=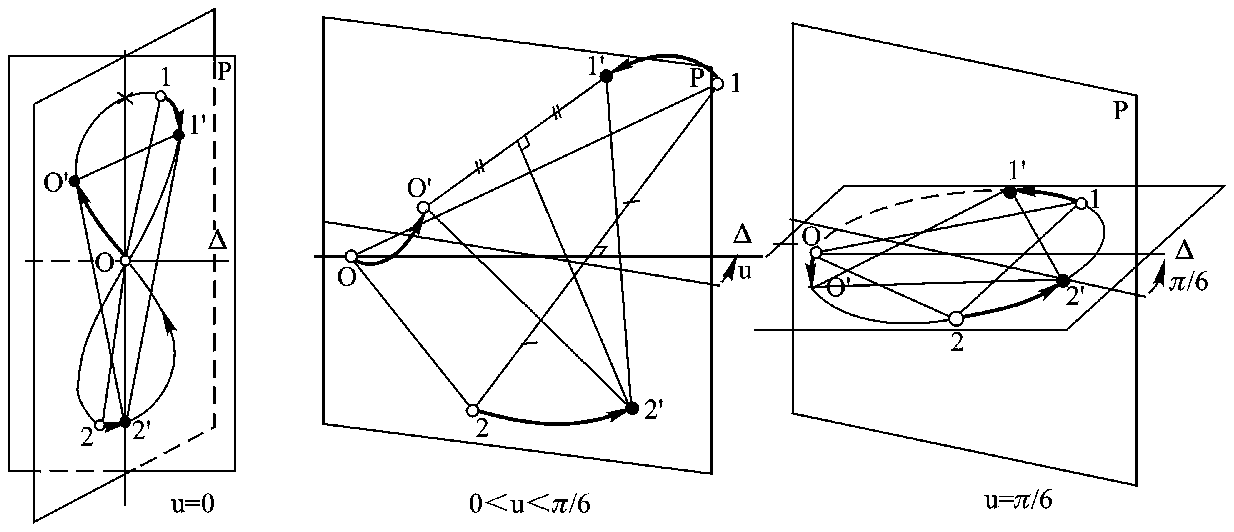,width=12cm}}}
\centerline{Figure 1 (fixed frame)}
\endinsert

For $u=0$, one gets the Eight in the vertical plane orthogonal to
$\Delta$ (and hence to $P$); for $u={\pi\over 6}$, one gets Lagrange
solution in the horizontal plane (containing $\Delta$ and orthogonal to
$P$). In a frame rotating around the vertical axis of an angle $-u$ in
time $T/12$, one gets a family of $D_6$-symmetric choreographies of
period $T$ between the Eight and twice Lagrange (figure~2). \midinsert
\centerline{\hbox{\psfig{figure=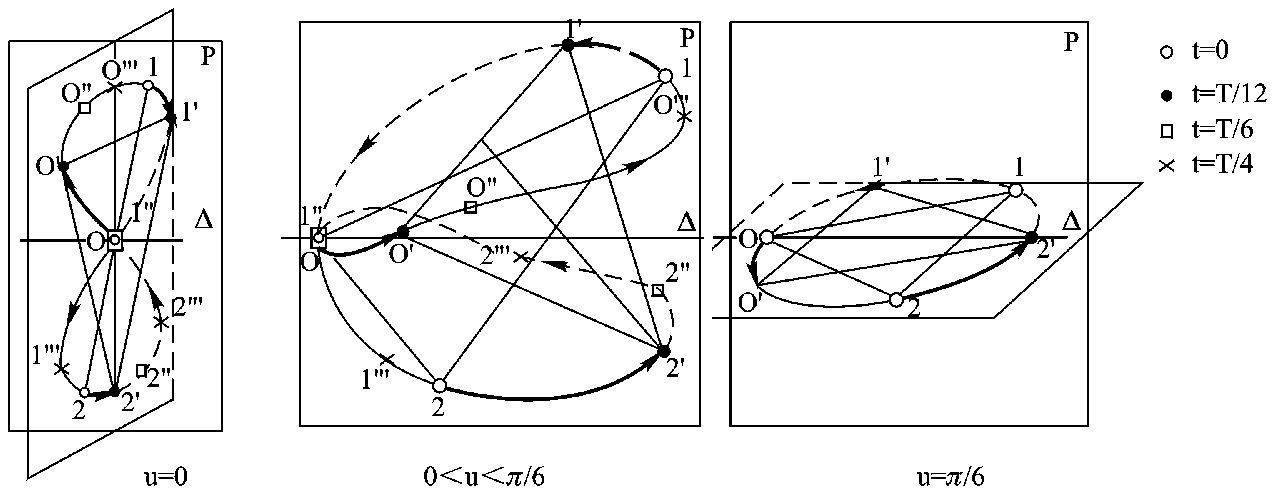,width=12cm}}}

\centerline{Figure 2 (rotating rame)}
\endinsert

 The relevant action of $D_6$ on the configuration space of three bodies in $\R^3$ is a direct
generalization of the one which leaves the Eight invariant. It is defined as follows (the notations
are the ones of [C3]):
$$\eqalign{&\alpha(s)(\vec r_0,\vec r_1,\vec r_2)=(\Sigma\vec r_2,\Sigma\vec r_0,\Sigma\vec r_1),
\quad \beta(s)(t)=t+T/6,\cr
&\alpha(\sigma)(\vec r_0,\vec r_1,\vec r_2)=(\Delta\vec r_0,\Delta\vec r_2,\Delta\vec r_1),\quad
\beta(\sigma)(t)=-t,}$$
where $\Sigma$ (resp. $\Delta$) denotes the symmetry with respect to the horizontal plane (resp. to the line
$\Delta$).

Thanks to the fact that a minimizer of the fixed-ends problem has no
collision, we need only show that a minimizing path has no collisions at
its ends, which can be done using local deformations as in one of the
proofs of the test assertion for the Kepler problem. The surprise is
that, using as a model the horizontal Lagrange family (which satisfies
the symmetry requirements), one can give a simple direct proof of the
absence of collisions in a minimizer:

1) the action of an admissible path undergoing a collision is bigger
than the action $\hat A_2=2^{-{5\over 3}}3^{2\over 3}\pi^{2\over
3}T^{1\over 3}$ (masses =1) of the horizontal relative equilibrium
solution $x_0$ of an equilateral triangle which rotates by $\pi\over 3$
in the same amount of time $T/12$;

 2) this last action is, for any $u\le{\pi\over 3}$, bigger than the one
$A(u)=\hat A_2\bigl[{3\over \pi}({\pi\over 3}-u)\bigr]^{2\over 3}$ of
the horizontal relative equilibrium solution $x_u$ of an equilateral
triangle which rotates by an angle $({\pi\over 3}-u)$ during the given
amount of time.

The first estimation, better than the one in [CM] ($\hat A_2=2^{1\over
3}A_2$) , appears in [ZZ2] in the case of the Eight. It follows from the
remark (at the basis of [V1] and [ZZ1]) that the action of a 3-body
problem splits into the sum of three terms, each of which is one third
of the action of the Kepler problem with attraction constant equal to
the total mass $M=3$. As the configurations at $t$ and $t+T/2$ are
symmetric with respect to the horizontal plane (compute $\alpha(s^3)$
and $\beta(s^3)$), any collision which occurs at $t_0$ occurs also at
$t_0+T/2$. The lower bound of the Kepler action during a period $T$ is
then twice the minimum of the Kepler action of an ejection-collision wih
attraction constant 3 and period $T/2$. But this is exactly the action
of $x_0$ during the period.

Finally, we prove that, for $0\le u<{\pi\over 6}$, the Lagrange solution
$x_u$ is not a minimizer. This is because the value $d^2{\cal
A}(x_u)(\xi,\xi)$ of the Hessian of the action on the vertical variation
$$\xi=\left(\sin ({2\pi t\over T}),\; \sin ({2\pi t\over T}+{2\pi\over 3}),\;
\sin ({2\pi t\over T}+{4\pi\over 3})\right)$$ which ``opens'' $x_u$ in
the direction of the Eight, is negative for $u<\pi/6$ and positive for
$u>\pi/6$. Indeed, the Hessian of $x_u$ is positive when $\pi/6<u\leq
\pi/3$, which supports Marchal's claim that $x_u$ is the minimizer when
$\pi/6\leq u\leq \pi/3$ (notice that its size increases to infinity and
its action decreases to 0 when $u$ tends to $\pi/3$).

\noindent {\bf Questions.} 1) Prove that for $u=0$, (the) minimum is
planar, hence (the) Eight.

\noindent 2) Our argument works for one value of $u$ at a time. As no
uniqueness is proved, neither is continuity with $u$ of the family. Such
continuity would imply the existence among the family of spatial 3-body
choreographies in the fixed frame.

\noindent{\bf Remark.} The first continuation of the Eight into a family
of rotating planar choreographies was given by Michel H\'enon [CGMS]
using the same program as in [H]. A third family should exist, rotating
around an axis orthogonal to the first two.

\head{5. Related results and open problems}

Two global questions seem to be out of reach at the moment:
unicity and possible extra symmetries of minimizers.

As an example of the first, numerical evidence by Sim\'o suggests
unicity of the Eight but in [CM] we do not even prove that each lobe is
convex, only that it is star-shaped (the problem is near the crossing
point). This is nevertheless enough to imply that the braid it defines
in space time $\R^2\times \R/T\Z$ (equivalent to the homotopy class in
the configuration space) is the ``Borromean rings", the signature of a
truly triple interaction (also noticed in [Ber] in a different context).

 As an example of the second one, we do not know if the $\Z/4\Z$-symmetry and the ``brake" property of the
Hip-Hop solution [CV] follow automatically from minimizing the action among loops such that $x(t+T/2)=-x(t)$
(compare 4.2 (ii)). One is tempted to compare this problem to the celebrated result of Alain Albouy [A2]
which states the existence of some symmetry in any central configuration of 4 equal masses (and implies that
there is only a finite number of them). But there is Moeckel's numerical example
[Mo2]  of a central configuration of eight equal masses without any symmetry. And according to [SW], there exists such an example
minimizing $\tilde U$ for $n=46$. For more on symmetry, see [V2].
\smallskip
Identifying minimizers, even when one knows that they are collision-free, is usually too difficult a problem
(see 4.2 (i) and (ii)). Understanding their stability properties may sometimes be
attempted theoretically [Ar],[O], or numerically [S1].
\smallskip\goodbreak
Another type of questions is connected with minimization with mixed constraints: symmetry
and homology or homotopy. One can ask, for example, if the Eight is a minimizing
choreography in its homology class $(0,0,0)$ (each side of the triangle has
zero total rotation). An interesting example of mixed
conditions may be found in [V2] where generalizations of the Hip-Hop lead to spatial choreographies of 4
equal masses.  But, as for most choreographies, no proof was found of the existence of Gerver's ``supereight"
with four equal masses [CGMS], [C3].
\smallskip
{\it I am indebted to Christian Marchal, Richard Montgomery, David Sauzin, Susanna Terracini and Andrea Venturelli for many illuminating
discussions and comments.}

\head{References}

\parindent 15mm

\item{[A1]} Albouy A. Lectures on the two-body problem, {\it Classical and Celestial Mechanics: The Recife
Lectures, H. Cabral F. Diacu ed.} (in press at Princeton University Press).

\item{[A2]} Albouy A. Sym\'etrie des configurations centrales de quatre
corps, {\it C. R. Acad. Sci. Paris, 320, 217--220 (1995)} \& The
symmetric central configurations of four equal masses, {\it
Contemporary Mathematics, vol. 198} (1996), 131--135.

\item{[AC]} Albouy A. and Chenciner A., Le
probl\` eme des $n$ corps et les distances mutuelles, {\it
Inventiones Mathematic\ae, 131}  (1998), 151--184.

\item{[AGT]} Arioli G., Gazzola F. and Terracini S. Minimization properties
of Hill's orbits and applications to some $N$-body problems, {\it
preprint}, (October 1999).

\item{[Ar]} Arnaud M.C.  On the type of certain periodic orbits minimizing the Lagrangian action
{\it Nonlinearity 11} (1998), 143--150.

\item{[B]} Betti E. Teorica delle forze newtoniane e sue applicazioni all'
elestrotatica e al magnetismo, {\it Pisa, Nistri} (1879).

\item{[Ber]} Berger M.A. Hamiltonian dynamics generated by Vassiliev
invariants,  {\it Journal of Physics A: Math. Gen. 34} (2001),
1363--1374.

\item{[Ch]} Chen K.C.  Action minimizing orbits in the parallelogram
four-body problem with equal masses, {\it Arch. Ration. Mech.
Anal., 158, no. 4} (2001), 293--318.

\item{[C1]} Chenciner A. Introduction to the N-body problem, {\it
Ravello summer school} (09-1997), http://www.bdl.fr/Equipes/ASD/person/chenciner/chenciner. htm

\item{[C2]} Chenciner A. Collisions totales, Mouvements compl\`etement
paraboliques et r\'eduction des homoth\'eties dans le probl\`eme
des $n$ corps, {\it Regular and chaotic dynamics, V.3, 3} (1998),
93--106.

\item{[C3]} Chenciner A. Action minimizing periodic solutions of the
$n$-body problem,  {\it ``Celestial Mechanics, dedicated to
Donald Saari for his 60th Birthday", A. Chenciner, R. Cushman, C.
Robinson, Z.J. Xia ed., Contemporary Mathematics 292}  (2002),
71--90.

\item{[C4]} Chenciner A. Simple non-planar periodic solutions of the $n$-body problem {\it Proceedings of the NDDS
Conference, Kyoto}, (2002).

\item{[CD]} Chenciner A. and Desolneux N. Minima de l'int\'egrale d'action
et \'equilibres relatifs de $n$ corps, {\it C.R. Acad. Sci.
Paris. t. 326, S\'erie I} (1998), 1209--1212. Corrections in {\it
C.R. Acad. Sci. Paris. t. 327, S\'erie I} (1998), 193 and in [C3].

\item{[CGMS]} Chenciner A., Gerver J., Montgomery R. and Sim\'o C. Simple
choreographies of $N$ bodies: a preliminary study, {\it Geometry,
Mechanics and Dynamics, Springer}, (to appear).

\item{[CM]} Chenciner A. and Montgomery R. A remarkable periodic solution of
the three body problem in the case of equal masses, {\it Annals
of Math., 152} (2000), 881--901.

\item{[CV]} Chenciner A. and Venturelli A. Minima de l'int\'egrale d'action
du Probl\`eme newtonien de $4$ corps de masses \'egales dans
$\R^3$ : orbites ``hip-hop", {\it Celestial Mechanics, vol. 77}
(2000), 139--152.

\item{[C-Z]} Coti Zelati V. Periodic solutions for $N$-body type problems,
{\it Ann. Inst. H. Poincar\'e, Anal. Non Lin\'eaire, v. 7, no. 5}
(1990), 477--492.

\item{[DGM]} Degiovanni M., Giannoni F. and Marino A., Periodic solutions of dynamical
systems with Newtonian type potentials,  {\it Ann. Scuola Norm.
Sup. Pisa Cl. Sci. 15} (1988), 467--494.

\item{[G]} Gordon W.B. A Minimizing Property of Keplerian Orbits, {\it American
Journal of Math. vol. 99, no. 15} (1977), 961--971.

\item{[H]} H\'enon M.  Families of periodic orbits in the three-body problem, {\it
Celestial Mechanics 10} (1974), 375--388.

\item{[M1]} Marchal C. The family $P_{12}$ of the three-body problem. The
simplest family of periodic orbits with twelve symmetries per
period, {\it Fifth Alexander von Humboldt Colloquium for
Celestial Mechanics}, (2000).

\item{[M2]} Marchal C. How the method of minimization of action avoids
singularities,  {\it Celestial Mechanics and Dynamical
Astronomy}, (to appear).

\item{[M3]} Marchal C. Handwritten supplement to the above paper and private
discussions.

\item{[Mo1]} Moeckel R. On central configurations, {\it Math. Z. 205}
(1990), 499--517.

\item{[Mo2]} Moeckel R. Some Relative Equilibria of N Equal Masses,
N=4,5,6,7; Addendum: N=8; {\it unpublished paper describing
numerical experiments} ($\leq 1990$).

\item{[M]} Montgomery R. Action spectrum and collisions in the three-body
problem,  {\it  ``Celestial Mechanics, dedicated to Donald Saari
for his 60th Birthday", A. Chenciner, R. Cushman, C. Robinson,
Z.J. Xia ed., Contemporary Mathematics 292} (2002), 173--184.

\item{[O]} Offin D. {\it Maslov index and instability of periodic orbits in
Hamiltonian systems}, preprint, 2002.

\item{[P]} Poincar\'e H. Sur les solutions p\'eriodiques et le principe de
moindre action,  {\it C.R.A.S. t. 123} (1896), 915--918.

\item{[SeT]} Serra E. and Terracini S. Collisionless Periodic Solutions to
Some Three-Body Problems, {\it Arch. Rational Mech. Anal., 120}
(1992), 305--325.

\item{[S1]} Sim\'o C. Dynamical properties of the figure eight solution of
the three-body problem, {\it ``Celestial Mechanics, dedicated to
Donald Saari for his 60th Birthday", A. Chenciner, R. Cushman, C.
Robinson, Z.J. Xia ed., Contemporary Mathematics 292},  (2002),
209--228.

\item{[S2]} Sim\'o C. New families of Solutions in $N$-Body Problems,
{\it Proceedings of the Third European Congress of Mathematics,
C. Casacuberta et al. eds. Progress in Mathematics, 201} (2001),
101--115.

\item{[SW]} Slaminka E. \& Woerner K. Central configurations and a theorem of Palmore {\it
Celestial Mechanics 48} (1990), 347--355.

\item{[V1]} Venturelli A. Une caract\'erisation variationnelle des
solutions de Lagrange du pro\-bl\`eme plan des trois corps, {\it
C.R. Acad. Sci. Paris, t. 332, S\'erie I} (2001), 641--644.

\item{[V2]} Venturelli A. Th\`ese, {\it Paris} (to be defended in
2002).

\item{[ZZ1]} Zhang S. \& Zhou Q., A Minimizing Property of
Lagrangian Solutions, {\it Acta Mathematica Sinica, English Series},
Vol. 17, No.3 (2001), 497--500.

\item{[ZZ2]} Zhang S. \& Zhou Q., Variational method for the
choreography solution to the three-body problem, {\it Science in China
(Series A)}, Vol. 45 No. 5 (2002), 594--597.

\vfill\eject
\end